\documentclass[graybox]{svmult}

\usepackage{mathptmx}       
\usepackage{helvet}         
\usepackage{courier}        
\usepackage{type1cm}        
\usepackage{makeidx}         
\usepackage{graphicx}        
\usepackage{multicol}        
\usepackage[bottom]{footmisc}

\usepackage{cite}
\usepackage{fancyvrb}        

\usepackage{amsmath}
\usepackage{amssymb}

\usepackage[european]{circuitikz}
\usepackage{relsize}

\makeindex             

\newcommand{\state}{x}

\newcommand{\statedisturbed}{\overline{\state}}

\newcommand{\control}{u}
\newcommand{\controlstar}{\control^\star}
\newcommand{\controldisturbed}{\overline{\control}}

\newcommand{\disturbancezero}{w}

\newcommand{\disturbance}{\overline{\disturbancezero}}

\newcommand{\costfunction}{J}
\newcommand{\costfunctionN}{\costfunction_N}

\newcommand{\openloopindex}{k}

\newcommand{\controlset}{U}

\newcommand{\controlsetconstrained}{\mathbb{\controlset}}

\begin{document}

\title*{Performance of Sensitivity based NMPC Updates in Automotive Applications}
\author{J\"{u}rgen Pannek and Matthias Gerdts}
\institute{J\"{u}rgen Pannek \at University of the Federal Armed Forces Munich, Werner-Heisenberg-Weg 39, 85577 Neubiberg/M\"{u}nchen, \email{juergen.pannek@unibw.de}
\and Matthias Gerdts \at University of the Federal Armed Forces Munich, Werner-Heisenberg-Weg 39, 85577 Neubiberg/M\"{u}nchen, \email{matthias.gerdts@unibw.de}}
%
%
\maketitle

\abstract*{In this work we consider a half car model which is subject to unknown but measurable disturbances. To control this system, we impose a combination of model predictive control without stabilizing terminal constraints or cost to generate a nominal solution and sensitivity updates to handle the disturbances. For this approach, stability of the resulting closed loop can be guaranteed using a relaxed Lyapunov argument on the nominal system and Lipschitz conditions on the open loop change of the optimal value function and the stage costs. For the considered example, the proposed approach is realtime applicable and corresponding results show significant performance improvements of the updated solution with respect to comfort and handling properties.}

\abstract{In this work we consider a half car model which is subject to unknown but measurable disturbances. To control this system, we impose a combination of model predictive control without stabilizing terminal constraints or cost to generate a nominal solution and sensitivity updates to handle the disturbances. For this approach, stability of the resulting closed loop can be guaranteed using a relaxed Lyapunov argument on the nominal system and Lipschitz conditions on the open loop change of the optimal value function and the stage costs. For the considered example, the proposed approach is realtime applicable and corresponding results show significant performance improvements of the updated solution with respect to comfort and handling properties.}

\section{Introduction}
\label{p41:sec_introduction}

Within the last decades, model predictive control (MPC) has grown mature for both linear and nonlinear systems, see, e.g., \cite{CamachoBordons2004, RawlingsMayne2009, GruenePannek2011}. 
Although analytically and numerically challenging, the method itself is attractive due to its simplicity and approximates an infinite horizon optimal control as follows:
In a first step, a measurement of the current system state is obtained which in the second step is used to compute an optimal control over a finite optimization horizon. 
In the third and last step, a portion of this control is applied to the process and the entire problem is shifted forward in time rendering the scheme to be iteratively applicable.

Unfortunately, stability and optimality of the closed loop may be lost due to considering finite horizons only. To ensure stability of the resulting closed loop, one may impose terminal point constraints as shown in \cite{KeerthiGilbert1988, Alamir2006} or Lyapunov type terminal costs and terminal regions, see \cite{ChenAllgoewer1998, MayneRawlingsRaoScokaert2000}. A third approach uses a relaxed Lyapunov condition presented in \cite{GrueneRantzer2008} which can be shown to hold if the system is controllable in terms of the stage costs \cite{Gruene2009, GruenePannekSeehaferWorthmann2010}. Additionally, this method allows for computing an estimate on the degree of suboptimality with respect to the infinite horizon controller, see also \cite{ShammaXiong1997, NevisticPrimbs1997} for earlier works on this topic.

Here, we use an extension of the third approach to the case of parametric control systems and subsequent disturbance rejection updates. In particular, we focus on updating the MPC control law via sensitivities introduced in \cite{Fiacco1983}. Such updates have been analysed extensively for the case of open loop optimal controls, see, e.g, \cite{GroetschelKrumkeRambau2001}, but were also applied in the MPC closed loop context in \cite{ZavalaBiegler2009, PannekGerdts2012}. In order to avoid the usage of stabilizing Lyapunov type terminal costs and terminal regions and obtain performance results with respect to the infinite horizon controller, we utilize results from \cite{PannekGerdts2012} in an advanced step setting, see, e.g., \cite{FindeisenAllgoewer2004}. 

In the following, we present the considered half car model from \cite{SpeckertDresslerRuf2009, PoppSchiehlen2010} and the imposed MPC setup. The obtained numerical results show that this approach is both realtime applicable and provides a cheap and yet significant performance improvement with respect to the comfort and handling objectives requested by our industrial partners.

\section{Problem setting}
\label{p41:sec_setting}

\begin{figure}[!ht]
	\sidecaption[t]
	\includegraphics[width=0.62\textwidth]{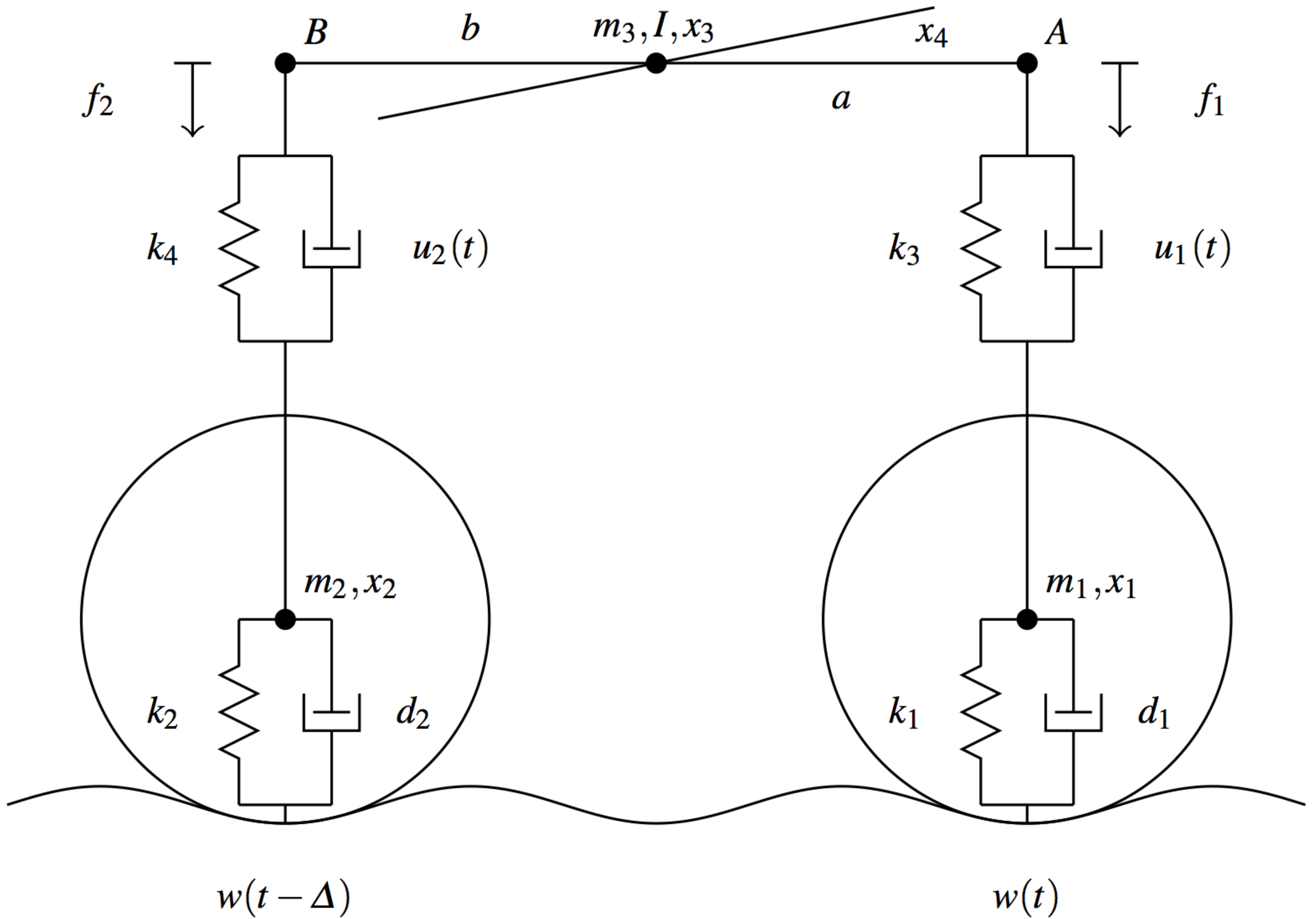}
		\caption{Schematical sketch of a halfcar subject to road excitation $\disturbancezero$}
		\label{fig:halfcar}
\end{figure}

Throughout this work we consider the control systems dynamics of a half car which originate from \cite{SpeckertDresslerRuf2009, PoppSchiehlen2010} and are slightly modified to incorporate active dampers, see Fig. \ref{fig:halfcar} for a schematical sketch. 
The resulting second order dynamics read
\begin{align}
	m_1 \ddot{\state}_1 & = m_1 g + f_3 - f_1 \qquad & m_3 \ddot{\state}_3 & = m_3 g - f_3 - f_4 \nonumber \\
	\label{dynamic}
	m_2 \ddot{\state}_2 & = m_2 g + f_4 - f_2 \qquad & I \ddot{\state}_4 & = \cos(\state_4) ( b f_3 - a f_4 )
\end{align}
where the control enters the forces
\begin{align*}
	f_1 & = k_1 ( \state_1 - \disturbancezero_1 ) + d_1 (\dot{\state}_1 - \dot{\disturbancezero}_1) \\
	f_2 & = k_2 ( \state_2 - \disturbancezero_2 ) + d_2 (\dot{\state}_2 - \dot{\disturbancezero}_2)  \\
	f_3 & = k_3 ( \state_3 - \state_1 - b \sin(\state_4) ) + \control_1 ( \dot{\state}_3 - \dot{\state}_1 - b \dot{\state}_4 \cos(\state_4) ) \\
	f_4 & = k_4 ( \state_3 - \state_2 + a \sin(\state_4) ) + \control_2 ( \dot{\state}_3 - \dot{\state}_2 + a \dot{\state}_4 \cos(\state_4) )
\end{align*}

Here, $\state_1$ and $\state_2$ denote the centers of gravity of the wheels, $\state_3$ the respective center of the chassis and $\state_4$ the pitch angle of the car. The disturbances $\disturbancezero_1$, $\disturbancezero_2$ are connected via $\disturbancezero_1(t) = \disturbancezero(t)$, $\disturbancezero_2(t) = \disturbancezero(t - \Delta)$ and the control constraints $\controlsetconstrained = [ 0.2 kNs/m, 5 kNs/m ]^2$ limit the range of the active dampers. The remaining constants of the halfcar are displayed in Table \ref{tab:halfcar}.
\begin{table}[!ht]
	\caption{Parameters for the halfcar example}
	\label{tab:halfcar}
	\begin{center}
		\begin{tabular}{cccc} \hline
			\hline\noalign{\smallskip}
			name & symbol & quantity & unit\\
			\noalign{\smallskip}\svhline\noalign{\smallskip}
			distance to joint & $a, b$ & $1$ & $m$ \\
			mass wheel & $m_1, m_2$ & $15$ & $kg$ \\
			mass chassis & $m_3$ & $750$ & $kg$\\
			inertia & $I$ & $500$ & $kg\, m^2$\\
			spring constant wheels & $k_1, k_2$ & $2\cdot 10^5$ & $kN/m$\\
			damper constant wheels & $d_1, d_2$ & $2\cdot 10^2$ & $kNs/m$\\
			spring constant chassis & $k_3, k_4$ & $1 \cdot 10^5$ & $kN/m$\\
			gravitational constant & $g$ & $9.81$ & $m/s^2$ \\
			\noalign{\smallskip}\hline
		\end{tabular}
	\end{center}
\end{table}

\section{MPC Algorithm}
\label{p41:sec_mpc}

In order to design a feedback for the half car problem \eqref{dynamic}, we impose the cost functional
\begin{align}
	\label{costfunctional}
	\costfunctionN(\state, \control, \disturbancezero) := \sum_{\openloopindex = 0}^{N-1} \mu_R F_{R}( \openloopindex ) + \mu_A F_{A}( \openloopindex )
\end{align}
following ISO 2631 with horizon length $N = 5$. The handling objective is implemented via
\begin{align*}
	& F_{R}( \openloopindex ) := \sum_{i = 1}^2 \int\limits_{\openloopindex T}^{(\openloopindex+1) T} \left( \frac{[ k_i ( \state_i(t) - \disturbancezero_i(t) ) + d_i ( \dot{\state}_i(t) - \dot{\disturbancezero}_i(t) ) ] - F_i}{F_i}  \right)^2 \, dt
\end{align*}
with nominal forces 
\begin{align*}
	F_1 & = ( a \cdot g \cdot ( m_1 + m_2 + m_3 ) ) / ( a + b ) \\
	F_2 & = ( b \cdot g \cdot ( m_1 + m_2 + m_3 ) ) / ( a + b )
\end{align*}
whereas minimizing the chassis jerk
\begin{align*}
	F_{A}(\openloopindex) := \int\limits_{\openloopindex T}^{(\openloopindex+1)T} \left( m_3 \dddot{\state}_3(t) \right)^2 \, dt
\end{align*}
is used to treat the comfort objective. Both integrals are equally weighted via $\mu_R = \mu_A = 1$ and are evaluated using a constant sampling rate of $T = 0.1s$ during which the control are held constant, i.e. the control is implemented in a zero--order hold manner. The nominal disturbance $\disturbancezero(\cdot)$ and the corresponding derivates are computed from road profile measurements taken at a sampling rate of $0.002s$ via a fast Fourier transformation (FFT). 

For the resulting finite time optimal control problem, we denote a minimizer of \eqref{costfunctional} satisfying all constraints by $\controlstar(\cdot, \state, \disturbancezero)$. Since the control must be readily computed at the time instant it is supposed to be applied, $\controlstar(\cdot, \state, \disturbancezero)$ is computed in an advanced step setting, cf. \cite{FindeisenAllgoewer2004}. To this end, the initial state $\state$ of the optimal control problem is predicted for a future time instant using the last known measurement and the intermediate control which is readily available from previous MPC iteration steps.

Since we want to apply sensitivity updates in case of measurement/prediction deviations and disturbances, we additionally precompute sensitivity information along the optimal open loop solution with respect to the predicted state $\partial \controlstar/\partial \state(\cdot, \state, \disturbancezero)$ and the nominal disturbances $\partial \controlstar/\partial \disturbancezero(\cdot, \state, \disturbancezero)$. Then, once the nominal control $\controlstar(\cdot, \state, \disturbancezero)$ is to be applied, we use newly obtained state and disturbance information $\statedisturbed$, $\disturbance$ to update the control via
\begin{align}
	\label{update formula}
	\controldisturbed(\cdot, \statedisturbed, \disturbance) := \controlstar(\cdot, \state, \disturbancezero) + 
	\begin{pmatrix}
		\frac{\partial \controlstar}{\partial \state}(\cdot, \state, \disturbancezero) \\
		\frac{\partial \controlstar}{\partial \disturbancezero}(\cdot, \state, \disturbancezero)
	\end{pmatrix}^\top
	\begin{pmatrix}
		\statedisturbed(\cdot) - \state(\cdot) \\
		\disturbance(\cdot) - \disturbancezero(\cdot)
	\end{pmatrix},
\end{align}
see also \cite{GroetschelKrumkeRambau2001, Fiacco1983} for details on the computation and limitations of sensitivities. 

For simplicity of exposition, we predict the initial state $\state$ using two sampling intervals $T$ of the closed loop control. Note that although larger predictions are possible, robustness problems are more likely to occur since predicted and real solutions usually diverge, see, e.g., \cite{LimonAlamoRaimondoBravoMunoyFerramoscaCamacho2009, FindeisenGruenePannekVarutti2011}.

\section{Numerical Results}
\label{p41:sec_numericalresults}

\begin{figure}[!ht]
	\sidecaption[t]
	\includegraphics[width=0.62\textwidth]{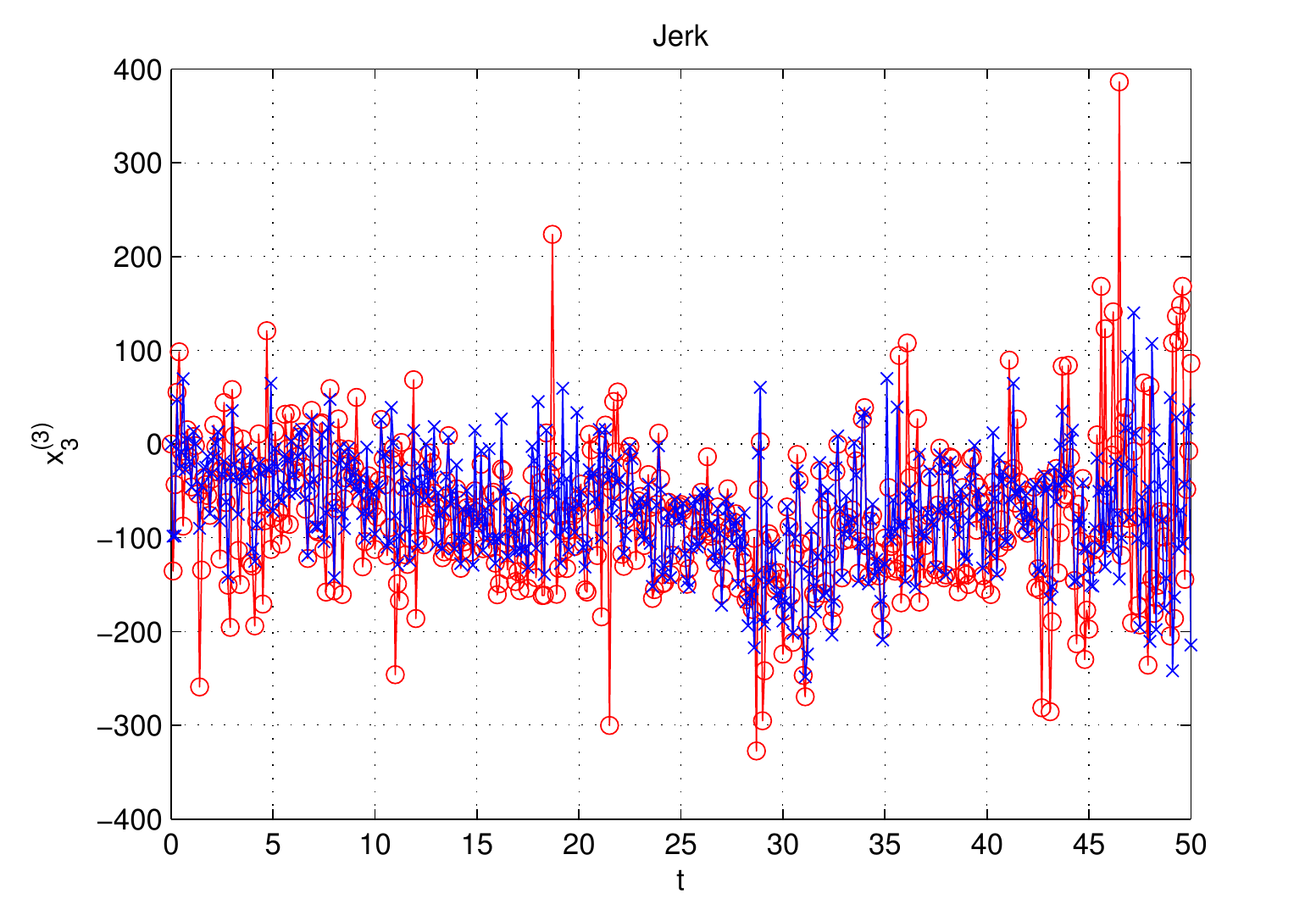}
	\caption{Comparison plot for the chassis jerk using MPC with ($\times$) and without sensitivity update ($\circ$).}
	\label{fig:results}
\end{figure}

During our simulations, we modified both the states of the system and the road profile measurements using a disturbance which is uniformly distributed in the interval $[-0.025 m, 0.025 m]$. For this setting, precomputation of $\controlstar(\cdot, \state, \disturbancezero)$, $\partial \controlstar/\partial \state(\cdot, \state, \disturbancezero)$ and $\partial \controlstar/\partial \disturbancezero(\cdot, \state, \disturbancezero)$ required at maximum $0.168s < 2 T = 0.2s$ which renders the scheme realtime applicable. As expected, the updated control law shows a better performance than the nominal control. The improvement cannot only be observed from Fig. \ref{fig:results}, but also in terms of the closed loop costs: For the considered race track road data we obtain an improvement of approximately $8.2 \%$ using the sensitivity update \eqref{update formula}. Although this seems to be a fairly small improvement, the best possible result obtained by a full reoptimization reveals a reduction of approximately $10.5 \%$ of the closed loop costs.

Note that due to the presence of constraints it is a priori unknown whether the conditions of the Sensitivity Theorem of \cite{Fiacco1983} hold at each visited point along the closed loop. Such an occurrance can be detected online by checking for violations of constraints or changes in the control structure. Yet, due to the structure of the MPC algorithm, such an event has to be treated if one of the constraints is violated at open loop time instant $\openloopindex = 1$ only which was not the case for our example.

\begin{acknowledgement}
This work was partially funded by the German Federal Ministry of Education and Research (BMBF), grant no. 03MS633G. 
\end{acknowledgement}

%
%
%



\end{document}